\input amstex
\documentstyle{amsppt}
\magnification=1200
\vsize19.5cm
\hsize13.5cm
\TagsOnRight
\pageno=1
\baselineskip=15.0pt
\parskip=3pt

\def\p{\partial}
\def\noo{\noindent}

\def\lam{\lambda}
\def\Om{\Omega}

\def\R{\bold R}

\def\wtt{\tilde}
\def\back{\backslash}
\def\Ga{\Gamma}

\def\dist{\text{dist}}

\def\ol{\overline}

\def\D{\nabla}
\def\phi{\varphi}

\def\M{\Cal M}

\nologo
 \NoRunningHeads

\topmatter

\title{A priori estimates and existence for\\
 a class of fully nonlinear elliptic equations \\
  in conformal geometry}\endtitle

\author{Xu-Jia Wang }\endauthor

\affil{The Australian National University }\endaffil

\address
 Centre for Mathematics and its Applications, Australian National
    University, \newline Canberra ACT 0200, Australia
 \endaddress
\email wang\@maths.anu.edu.au \endemail

\thanks{
This work was supported by the Australian Research Council.
\newline
$\text{\ \ \ \ \,}$ Published in {\it Chinese Annals of
Mathematics}, 27(B) (2006), 169-178.
 }\endthanks

\abstract {In this paper we prove the interior gradient and second
derivative estimates for a class of fully nonlinear elliptic
equations determined by symmetric functions of eigenvalues of the
Ricci or Schouten tensors. As an application we prove the
existence of solutions to the equations when the manifold is
locally conformally flat or the Ricci curvature is positive.
}\endabstract


\endtopmatter

\vskip-10pt

\document

\baselineskip=13.2pt
\parskip=3pt

\centerline{\it Dedicated to the memory of Professor S.S. Chern}

 \vskip20pt

\centerline {\bf 1. Introduction}

\vskip10pt

Let $(\M, g_0)$ be a compact Riemannian manifold of dimension
$n\ge 3$. Denote by $[g_0]$ the set of metrics conformal to $g_0$.
For $g\in [g_0]$ we denote by $R$, $Ric$, and $A_{g}=\frac
1{n-2}(Ric_{g}-\frac {R_{g}} {2(n-1)}g)$ the scalar curvature, the
Ricci tensor, the Schouten tensor, respectively. Let
$\lam(A_g)=(\lam_1, \cdots, \lam_n)$ denote the eigenvalues of
$A_g$ with respect to $g$. In this paper we study the regularity
and existence of solutions to the equation
$$f(\lam)=\phi,\tag 1.1$$
where $\phi$ is a positive, smooth function.

When $f=\sum \lam_i$, equation (1.1) is the prescribing scalar
curvature equation. When $\phi\equiv 1$, it is the Yamabe problem.
In this paper we assume that $f$ is a nonlinear function defined
on an open convex cone $\Ga\subset\R^n$ and satisfying the
following conditions:
\item {}(f1) $f>0$ in $\Ga$ and $f=0$ on $\p\Ga$;
\item {}(f2) $f_i>0$ in $\Ga$, where $f_i$ denotes derivative in
the $i$th variables;
\item {}(f3) $f$ is concave;
\item {}(f4) $f$ is invariant under exchange of variables;
\item {}(f5) $f$ is homogeneous of degree $\alpha$ for some $\alpha>0$,
namely $f(t\lam)=t^\alpha f(\lam)$ $\forall\ t>0$.

\noo By the concavity, we have $\alpha\in (0, 1]$. Note that (f1)
and (f5) also implies (f2). If $\phi\equiv 1$, then (f5) (with
$\alpha=1$) also follows from (f1)-(f4), as one can define a new
function $\wtt f$ such that $\{\lam\in\R^n\ |\ \wtt
f(\lam)=1\}=\{\lam\in\R^n\ |\ f(\lam)=1\}$ [LL2]. We assume the
cone $\Ga$ satisfies
\item {}(g1) $\Ga_n\subset\Ga\subset \Ga_1$, where
$\Ga_1$ and $\Ga_n$ are given in (1.3) below;
\item {}(g2) if $\lam\in\Ga$, then any permutation of $\lam$
also lies in $\Ga$.

Equation (1.1) is referred to as conformal Hessian equation. A
related equation is the Hessian equation, that is when $\lam$ are
eigenvalues of the Hessian matrix $\D^2 u$. The Hessian equation
was first studied in [CNS, I]; see also [CW, TW1] for further
results. There are many functions satisfy the above conditions,
for example,

\noo (i) The conformal $k$-Hessian equation. Let
$$\sigma_k(\lam)
  =\sum_{i_1<\cdots<i_k}\lam_{i_1}\cdots\lam_{i_k}  \tag 1.2$$
be the $k$-th elementary symmetric polynomial and let
$f(\lam)=\sigma_k^{1/k}(\lam)$. The corresponding cone is given by
$$\Ga_k=\{\lam\in \R^n\ |\
        \sigma_j(\lam)>0\ \text{for}\ j=1, \cdots, k\}.\tag 1.3$$
Note that $\Ga_1$ is the half-space $\{\lam\in\R^n\ |\
\sum\lam_i>0\}$ and $\Ga_n$ is the positive cone $\{\lam\in\R^n\
|\ \lam_i>0\ \ \forall\ 1\le i\le n\}$.

The quotient equation
$$f(\lam)=\big(\frac{\sigma_k}{\sigma_l}\big)^{1/(k-l)} (\lam)\tag 1.4$$
also satisfies (f1)-(f5), where $0\le l<k\le n$ and $\Ga=\Ga_k$.

\noo (ii) For any integer $1\le k\le n$ and $\delta\ge 0$,
$$f(\lam)=\delta\sum_{i=1}^n \lam_i+
  \min\{\lam_{i_1}+\cdots+\lam_{i_k}
   \ |\ 1\le i_1<\cdots<i_k\le n\} \tag 1.5$$
and $\Ga=\{\lam\in\R^n\ |\ f(\lam)>0\}$. When $\delta>0$, this is
indeed the Pucci operator, and is uniformly elliptic [GT].

\noo  (iii) More functions satisfy the above conditions, such as
$$f(\lam)=(\sum_{i=1}^n \lam_i^{-2})^{-1/2}, \tag 1.6$$
and
$$f(\lam)=\big[ \sum  \lam_1^{-\alpha_1}
    \cdots\lam_n^{-\alpha_n}\big]^{-1/k},  \tag 1.7 $$
where the sum is taken over all nonnegative integers $\alpha_1,
\cdots, \alpha_n$ with $\sum \alpha_i=k$, and $\Ga=\Ga_n$, see
[G]. One can find more functions satisfying the above conditions
[Kr].

\noo (iv) Let $f_1, f_2$ satisfy (f1)-(f5) and $f_2$ is
homogeneous of degree 1. Let
$$f(\lam)=f_1(\lam+\delta f_2(\lam) e),\tag 1.8$$
where $\delta>0$ is a constant, $e=(1, \cdots, 1)$.

From the definition of Schouten tensor, the Ricci curvature
$\mu=(\mu_1, \cdots, \mu_n)$ is given by
$$\mu_i=\lam_i+\frac{1}{n-2}\sum\lam_i \tag 1.9$$
Hence if we choose $\delta=\frac {1}{n-2}$ and $f_2(\lam)=\sum
\lam_i$ in (1.8), then we get an equation for the Ricci curvature
$\mu=(\mu_1, \cdots, \mu_n)$,
$$f_1(\mu)=1\tag 1.10$$

To write (1.1) as a partial differential equation, we use the
conformal transformation
$$g=v^{\frac{4}{n-2}}g_0. $$
Then (1.1) becomes
$$f(\lam(A^v))=\phi, \tag 1.11$$
where
$$A^v=\frac {2}{n-2}v^{-\frac{n+2}{n-2}}
  (\D^2_{conf} v+\frac {n-2}{2} vA_{g_0}).\tag 1.12$$
and
$$\D^2_{conf}v=-\D^2 v+\frac {n}{n-2}\frac {\D v\otimes\D v}{v}
  -\frac{1}{n-2}\frac {|\D v|^2}{v}I.\tag 1.13$$
In this paper we call $\D^2_{conf}v$ the {\it conformal Hessian
matrix}.

Equation (1.11) is a fully nonlinear partial differential
equation, in order that it is elliptic, we assume that the
eigenvalues $\lam\in\Ga$. Accordingly we say a metric $g$ {\it
$\Ga$-admissible}, or simply admissible, if $\lam\in\Ga$. Denote
by $[g_0]_\Ga$ the set of all conformal admissible metrics,
$$[g_0]_\Ga=\{g\in [g_0]\ |\ \lam(A_g)\in\Ga\}. \tag 1.14$$
When $g$ is an admissible solution to (1.1), the function $\phi$
must be positive. In the following we will consider solutions with
$\lam(A_g)\in \Ga$ only.

The existence of admissible solutions to (1.1) has been studied by
many authors and most papers are concentrated on the conformal
$k$-Hessian equation. Just to mention a few, the first existence
result is by [CGY], where the existence of solutions was proved
for the case $k=2, n=4$. Subsequently [LL1], [GW2] proved the
existence of solutions for general $k$ and $n$, but on locally
conformally flat manifolds. In [GV1] the authors proved the
existence for $k>\frac n2$, provided $\M$ is not conformal to the
unit sphere.

Recently with Sheng and Trudinger [STW], we proved the existence
of solutions to the conformal $k$-Hessian equation for $k\le \frac
n2$, provided equation (1.1) is variational. In [TW2] we proved
the compactness of the set $[g_0]_{\Ga_k}$ for $k>\frac n2$, for
any compact manifold not conformal to $S^n$, which also yields the
existence of solutions to the conformal $k$-Hessian equation for
$k>\frac n2$. Note that in [CGY, LL1, GW2, STW] the existence was
proved for $\phi\equiv 1$.

In this paper we establish the interior a priori estimates for
solutions to (1.1) under the mild conditions (f1)-(f5). The
interior a priori estimates were proved for the conformal
$k$-Hessian equation by P. Guan and G. Wang [GW1]. They also
proved the interior estimates for the quotient equation (1.4). Our
proof here uses a blow-up argument and is based on the Liouville
theorem of Y.Y. Li [L] and applies to the general equation (1.1)
with $f$ satisfying (f1)-(f5). From the a priori estimates, we
prove the existence of solutions of (1.1) when $\M$ is locally
conformally flat and $\phi\equiv 1$, or when the cone $\Ga$ is
smaller than $\Ga_{n/2}$. In the latter case, the Ricci curvature
must be positive.

This paper is arranged as follows. In Section 2 we state our main
results. In Section 3 we prove the interior second derivative
estimate. In Section 4 we prove the interior gradient estimate.
The final section 5 discusses that existence of solutions to
(1.1).

In a workshop at Berkeley in November 2005, the author learnt that
Sophie Chen has independently proved the interior a priori
estimates [Ch], see Remark 2.1 for more details. The author would
also like to thank Viaclovsky for discussions on Theorems 2.3 and
2.4 at the workshop.

\vskip30pt

\centerline{\bf 2. Main results}

\vskip10pt

For the a priori estimates in Theorems 2.1 and 2.2 below, we
assume that $\phi$ is a nonnegative function satisfying
$\D^2\phi\ge -C$ for some constant $C>0$.

\proclaim{Theorem 2.1} Let $v\in C^3$ be an admissible positive
solution of (1.11) in a geodesic ball $B_{r}(0)\subset\M$. Then
$$\frac {|\D v|}{v}(0)\le C,\tag 2.1$$
where $C$ depends only on $n$, $r$, $\phi$, $f$, $\sup v$, and
$g_0$, and $\D $ denotes the covariant derivative with respect to
$g_0$.
\endproclaim

\proclaim{Theorem 2.2} Let $v\in C^4$ be an admissible positive
solution of (1.11) in  $B_r(0)$. Then
$$|\D^2v|(0)\le C,\tag 2.2$$
where $C$ depends only on $n$, $r$, $\phi$, $f$, $\inf v$, $\sup
(v+|\D v|)$.
\endproclaim

\noo{\it Remark 2.1}. As indicated in the introduction, the author
learnt in a workshop at Berkeley in early November 2005
(31/10-4/11, 2005) that Sophie Chen [Ch] has independently proved
Theorems 2.1 and 2.2. Her preprint was available in September
2005, while the first version of the paper was finished in October
before the workshop. The proofs in [Ch] and in this paper are
completely different. The proof in this paper uses a blow-up
argument and the Liouville theorem [L], while the estimates in
[Ch] are derived directly from an auxiliary function and so are
more favourable.

Once the second derivatives are bounded, the equation becomes
uniformly elliptic if $\phi$ is positive, and further regularity
follows from the Evans-Krylov regularity theory.

The proof of the second derivative estimate is similar to that in
[GW1,GW3,STW]. For the interior gradient estimate, we use a
different approach. An important property of the conformal Hessian
equation is its invariance under Kelvin transformation (when the
equation is defined in $\R^n$). Hence by the moving plane method,
Y. Li [L] proved the Liouville theorem for the equation
$$F(\D^2_{conf} v)=0\ \ \ \text{in}\ \ \R^n .\tag 2.3$$
That is an entire positive solution of (2.3) must be a constant.
Using a blow-up argument, we prove the gradient estimate by the
Liouville theorem and the interior second derivative estimate.

As an application, we prove the existence of solutions to equation
(1.1) when $\M$ is locally conformal flat or the cone $\Ga$ is
smaller that $\Ga_{n/2}$.

\proclaim{Theorem 2.3} Let $(\M^n, g_0)$ be a compact, locally
conformally flat manifold. Suppose $\phi\equiv 1$. Then there
exists a solution to problem (1.1).
\endproclaim

\noo{\it Remark 2.2}. The existence of solutions to problem (1.1),
for general $f$ and $\phi\equiv 1$, was also proved in [LL2] (see
Theorem 1.1 in [LL2]).

We also remark that if $\phi$ is not a constant, there are
obstructions to the existence of solutions in general.

In our next existence result, we assume the cone $\Ga$ is relative
small, so that any admissible metric has positive Ricci curvature.
For this purpose we introduce a cone $\Sigma_\delta$ as in [GV2],
that is
$$\Sigma_\delta
  =\{\lam\in\R^n\ |\ \min \lam_i+\delta\sum\lam_i>0\},\tag 2.4$$
where $\delta\ge 0$. From (1.9),
$$Ric_g\ge 0\ \ \ \text{if and only if}\ \
          \lam(A_g)\in \Sigma_{\frac{1}{n-2}},\tag 2.5$$
and $Ric_g>0$ if and only if $\lam(A_g)\in \Sigma_\delta$ with
$\delta< \frac 1{n-2}$, except when $g$ is locally a flat metric.
Taking $l=1$ in the proof of Lemma 4.2 in [TW1], we also have
$$\Ga_k\subset \Sigma_{\frac {n-k}{n(k-1)}}.\tag 2.6$$
In particular we have $Ric_g\ge 0$ if $\lam(A_g)\in\Ga_{n/2}$ and
$Ric_g > 0$ if $\lam(A_g)\in\Ga_k$ for $k>\frac n2$.

\proclaim{Theorem 2.4} Let $(\M^n, g_0)$ be a compact manifold not
conformally equivalent to the unit sphere.  Suppose $\phi$ is a
positive smooth function and the cone $\Ga\subset \Sigma_\delta$
for some $0\le \delta<\frac 1{n-2}$. Then there exists a solution
to problem (1.1).
\endproclaim

\noo{\it Remark 2.3}. Theorem 2.4 was first proved by Gursky and
Viaclovsky [GV1] for the conformal $k$-Hessian equation for
$k>\frac n2$, see also [TW2] for a different proof. By the very
recent a priori estimates in [Ch], Gursky and Viaclovsky extended
their existence result to general curvature function $f$ as in
Theorem 2.4 [GV1].

\noo{\it Remark 2.4}. In Theorems 2.3 and 2.4, we will not only
prove the existence of solutions but also the compactness of the
set of solutions if $(\M, g_0)$ is not conformally equivalent to
the unit sphere. In particular under the assumptions in Theorem
2.4, the set of all admissible metrics (subject to the volume
constraint Vol$\M_g=1$) is compact.

\vskip30pt

\centerline{\bf 3. Second derivative estimate }

\vskip10pt

Denote
$$F(A^v)=f(\lam(A^v)), \tag 3.1$$
where $A^v$ is the matrix given in (1.12). Regard $F$ as a
function of $n\times n$ real symmetric matrices, $F=F(a_{ij})$.
Then from [CNS], we know that $F$ is elliptic by (f2), concave by
(f3), and symmetric under orthogonal transformations by (f4). That
is
$$ \{F^{ij}\}:=\{\frac {\p F}{\p a_{ij}}\}>0\tag 3.2$$
for any matrix $A=(a_{ij})$ with eigenvalues in $\Ga$,
$$ F^{ij, st}
 :=\frac {\p^2 F}{\p a_{ij}\p a_{st}}b_{ij}b_{st}\le 0\tag 3.3$$
for any matrix $B=(b_{ij})$, and
$$F(OAO')=F(A)\tag 3.4$$
for any orthogonal matrix $O$.

For the a priori estimates for the second derivatives, it is
convenient to use the conformal transformation $g=u^{-2}g_0$. Then
$$A_g=\frac {u_{ij}}{u}-\frac {|\D u|^2}{2u^2}g_0+ A_{g_0}$$
and equation (1.1) becomes
$$F(U)=\phi u^{-k\alpha}, \tag 3.5$$
where
$$U=\{u_{ij}-\frac {|\D u|^2}{2u}g_0+uA_{g_0}\}$$

Let $u\in C^{3.1}$ be an admissible positive solution of (3.5) in
a geodesic ball $B_{r}(0)\subset\M$. We want to prove
$$|\D^2u|(0)\le C,\tag 3.6$$
where $C$ depends only on $n$, $r$, $\sup(u+u^{-1})$, $\sup |\D
u|$, and $g_0$.

Let
$$z=\rho^2u_{\xi\xi}$$
where $\xi$ is any unit tangential vector,
$u_{\xi\xi}=u_{ij}\xi_i\xi_j$, $\rho(x)=(1-\frac {|x|^2}{r^2})^+$
is a cut-off function, $|x|$ denotes geodesic distance from $0$.
Assume that $\sup z$ is attained at $x_0$ and in direction $e=(1,
0, \cdots, 0)$. In an orthonormal frame at $x_0$, we may assume by
a rotation of axes that $\wtt u_{ij}:=u_{ij}+ua_{ij}$ is diagonal.
Then at $x_0$, $F^{ij}$ is diagonal and
$$\align
 0 & =(\log z)_i  =\frac {2\rho_i}\rho
  +\frac {\wtt u_{11,i}} {\wtt u_{11}},\tag 3.7\\
 0& \ge (\log z)_{ii} =(\frac {2\rho_{ii}}{\rho}
              -\frac{6\rho_i^2}{\rho^2})
  +\frac {\wtt u_{11, ii}}{\wtt u_{11}}.\tag 3.8\\
 \endalign   $$
Next, differentiating equation (3.5) twice,  we obtain
$$  F^{ij} U_{ij, kk}   =-\frac {\p^2 \mu(\sigma_k(\lam(U)))}
  {\p U_{ij}\p U_{rs}}U_{ij,k}U_{rs,k}+\D_k^2 (\phi u^{-k\alpha})
    \ge \D_k^2 (\phi u^{-k\alpha}),\tag 3.9$$
where $U_{ij, k}=\D_k U_{ij}$. Assume that
$$\frac{|\D u|}{u}\le C. \tag 3.10$$
Then by (3.7) and the Ricci identities,
$$U_{ij, 11}=u_{ij 11}-\frac {u_{k1}^2}{u}\delta_{ij}
                +O(\frac{1+u_{11}}{\rho})
 = u_{11 ij}-\frac {u_{k1}^2}{u}\delta_{ij}
                +O(\frac{1+u_{11}}{\rho}).\tag 3.11$$
Hence we obtain
$$\align
 0 & \ge \sum_i F^{ii}(\log z)_{ii}
  \ge -\frac {C}{\rho^2}\Cal F
 + F^{ii}\frac {\wtt u_{11,ii}}{\wtt u_{11}}\tag 3.12\\
   & \ge -\frac {C}{\rho^2}\Cal F
    +\frac {u_{11}^2}{u \wtt u_{11}}\Cal F
    +\frac {1}{\wtt u_{11}} \D_1^2 (\phi u^{-k\alpha}),\\
    \endalign  $$
where $\Cal F=\sum F^{ii}$.

Note that for any constant $a\ge 0$, there exists a constant
$C_a>0$ such that $\Cal F(A)>C_a$ when $F(A)=a$. To see this,
consider an arbitrary point $\lam_0\in\p\Ga$ satisfying
$\sup_{t>0} f(\lam_0+t\gamma)>a+1$, where $\gamma$ is the inner
unit normal of $\p\Ga$, or a proper unit vector at $\lam_0$,
pointing to the interior of $\Ga$ if $\p\Ga$ is not $C^1$ at
$\lam_0$. Let $t_a>0$ such that $f(\lam_0+t_a\gamma)=a$. Then the
derivative $\frac {d}{d t} f(\lam_0+t_a\gamma)$ is strictly
positive. Note that $\phi u^{-k\alpha}$ is bounded. Hence $\Cal
F\ge C$ in (3.12). Hence (3.6) holds.

As indicated before, the above proof is essentially contained in
[GW1,GW3,STW]. We include the proof for the convenience of the
readers.

We will establish the gradient estimate (3.10) in the next
section. Once the first and second order derivatives are bounded,
the equation becomes uniformly elliptic. Hence by Evans-Krylov's
regularity, we have the following regularity result.

\proclaim{Theorem 3.1} Let $u\in C^{3,1}$ be a positive solution
of (3.5) in $B_r(0)\subset \M$. Suppose $\phi>0, \in C^{1,1}$.
Then for any $\alpha\in (0, 1)$,
$$\|u\|_{C^{3, \alpha} (B_{r/2}(0))}\le C , \tag 3.13$$
where $C$ depends only on $n, r$, $\inf_\M u$, $\phi$, and $g_0$.
\endproclaim

\vskip30pt

\centerline{\bf 4. Interior gradient estimate}

\vskip10pt

In this section we prove an interior gradient estimate for
equation (1.11). We write equation (1.11) in the form
$$F(\D^2_{conf} v+A v)
     =\phi v^{\alpha\frac{n+2}{n-2}}, \tag 4.1$$
where $A=\frac{n-2}{2}A_{g_0}$. Let $v$ be a solution of (4.1) in
a geodesic ball $B_r(0)\subset\M$ and let
$$z=\rho \frac {|\D v|}{v}$$
be an auxiliary function, where $\rho=1-\frac {|x|^2}{r^2}$ is a
cut-off function. We want to prove that $z$ is uniformly bounded,
so that (2.1) holds

If estimate (2.1) is not true, there is a sequence of solutions
$v_k$ such that $\sup z_k\to\infty$, where $z_k=\rho \frac {|\D
v_k|}{v_k}$. Assume $\sup z_k$ is attained at $x_k$. We may assume
$$v_k(x_k)=1, $$
for otherwise we may replace $v_k$ by $\hat v_k=\frac
{v_k}{v_k(x_k)}$. Then $\hat v_k$ satisfies
$$F(\D^2_{conf}v +Av)=\ol c_k \phi v^{\alpha\frac{n+2}{n-2}} ,$$
where $\ol c_k=(v_k(x_k))^{4\alpha/(n-2)}$ and $\alpha>0$ is the
homogeneity constant in (f5). Note that our estimate allows $\sup
z$ depends on $\sup v$. Hence $\ol c_k$ is bounded above.

Denote $d_k=r-|x_k|=\dist (x_k, \p B_r(0))$. We have
$$d_k |\D v_k(x_k)|\ge r\sup z_k\to\infty.$$
Choose a normal coordinate at $x_k$ and identify $B_r(x_k)$ with a
Euclidean ball $B^e_r(0)$ by the exponential map. Make the
dilation
$$y=x|\D v_k(x_k)|.$$
Then $v_k$ is defined in a ball $B^e(0)$ of radius
$$r_k:=d_k|\D v_k(x_k)|\to\infty \tag 4.2$$
and satisfies the equation
$$F(\D^2_{conf}v +A_kv)=c_k\phi v^{\alpha\frac{n+2}{n-2}}, \tag 4.3$$
where $A_k=A |\D v_k(x_k)|^{-2}$ and $c_k=|\D
v_k(x_k)|^{-2\alpha}$. Moreover,
$$v_k=1\ \ \ |\D v_k|=1\ \ \ \text{at}\ \  y=0.\tag 4.4$$
Since $\sup z_k$ is attained at $y=0$ and note that $\rho>\frac 12
\rho(x_k)$ when $|x-x_k|<\frac 12 d_k$,  we see that $\frac {|\D_y
v_k|}{v_k}\le 2$ for $y\in B^e_{r_k/2}(0)$. Hence $v$ and $v^{-1}$
are locally uniformly bounded. Therefore for any $R>0$, by the
interior second derivative estimate, $v_k$ is uniformly bounded in
$C^{1,1}(B_R^e(0))$ provided $k$ is sufficiently large such that
$\frac {r_k}2>R$. Hence by the Arzela-Ascoli lemma, there is a
subsequence of $v_k(y)$ which converges to a limit function $v\in
C^{1,1}(\R^n)$, and $v$ is a solution of
$$F(\D^2_{conf}v  )=0. \tag 4.5 $$
Hence by [L], we have $v\equiv 1$. On the other hand, by (4.4) and
the interior second derivative estimate, we have $|\D v|=1$ at
$y=0$. We reach a contradiction and thus proved the interior
gradient estimate.

\vskip30pt

\centerline{\bf 5. Existence of solutions}

\vskip10pt

The existence of solutions to equation (1.1) for general $f$ has
been studied in [GV1,LL2] and [TW2], where the equation need not
to be variational. The key assumption is the interior a priori
estimates, which were established previously for the conformal
$k$-Hessian equations and their quotient equations (see [GW1]).
With our a priori estimates (Theorems 2.1 and 2.2) for general $f$
(see also [Ch]), we can prove the existence of solutions to (1.11)
as in [LL1, TW2], as stated in Theorems 2.3 and 2.4. So we outline
the proof here.

\vskip5pt

The proof of Theorem 2.3 is already in [LL1] (assuming the a
priori estimates, Theorems 2.1 and 2.2). Indeed, it suffices to
prove that the solution is uniformly bounded, namely
$$\sup v\le C\tag 5.1$$
for some $C$ depending only on $(\M, g_0)$ but independent of $f$.
We point out the main idea here. Let $\wtt\M$ be the universal
cover of $\M$, with the pull back of $g_0$ as the metric. Then the
function $\wtt v$ on $\wtt M$, determined by $v$, is also a
solution to (1.11). If $\wtt\M$ is the unit sphere, by the
stereographic projection, $\wtt\M$ (after taking away one point)
is conformally equivalent to the Euclidean space. Hence by the
Liouville Theorem [LL2], $\wtt v$ has a unique maximum point.
Hence $\sup v$ is bounded.

If $\wtt\M$ is not the unit sphere, it is conformal to a domain
$\wtt\Om\subset S^n$, $\ne S^n$ [SY]. Using the stereographic
projection, $\wtt\M$ is conformal to a domain $\Om\subset\R^n$,
$\ne R^n$. By the moving plane argument [Ye], the gradient of
$\wtt v$ is uniformly bounded, $|\D \wtt v|\le C$ for some $C$
depending only on $\Om$ (independent of $f$). Hence $v$ is
uniformly bounded on $\M$.

By (5.1) and Theorems 2.1 and 2.2, we obtain
$$\|v\|_{C^{2, \alpha}}
   +\|v^{-1}\|_{C^{2, \alpha}}\le C. \tag 5.2$$
Hence by a degree argument, we obtain a solution of (1.11).

\vskip5pt

For the proof of Theorem 2.4, we proceed as in [TW2]. We have a
stronger result, that is

\proclaim{Theorem 5.1} Let $(\M, g_0)$ be a compact $n$-manifold
not conformally equivalent to the unit sphere $S^n$. Suppose
$\Ga\subset\Sigma_\delta$ for some $\delta<\frac 1{n-2}$.  Then
the set $[g_0]_\Ga$ is compact in $C^0(\M)$ and satisfies the
following Harnack inequality, namely for any $g=e^{-2w} g_0\in
[g_0]_\Ga$,
$$|w(x)-w(y)|\le C|x-y|^\beta \tag 5.3$$
for some constant $C$ depending only on $\Ga$ and $(\M, g_0)$,
independent of $f$, where $\beta=\frac{1-\delta(n-2)}{1+\delta}$
and $|x-y|$ denotes geodesic distance in the metric $g_0$ between
$x$ and $y$.
\endproclaim

The compactness above is understood in the sense that $cg$ is
regarded as the same metric $g$ for any positive constant $c$.
Theorem 2.4 follows immediately from Theorem 5.1 by a degree
argument [TW2], and the proof of Theorem 5.1 is also similar to
that of Theorem A in [TW2]. We outline the main steps here.

\noo (i) For any metric $g=e^{-2 w} g_0\in [g_0]_\Ga$, by
subtracting a linear function we may assume $\sup w=0$. By the
assumption $\Ga\subset\Sigma_\delta$, one can show $u=e^w$ is
uniformly Holder continuous with Holder exponent $\beta
=\frac{1-\delta(n-2)}{1+\delta}$.

\noo (ii) If (5.3) is not true, there is a sequence of metrics
$g_k=e^{-2 w_k} g_0\in [g_0]_\Ga$ such that $\inf w_k\to -\infty$.
Let $w=\lim w_k$. By (i), $w$ is singular at some points. Let
$\wtt w(r)=\inf\{h\ |\ B_r(0)\subset \{w<h\}\}$ be the least
radial function satisfying $\wtt w\ge w$. We show that either
$\wtt w$ is H\"older continuous or
$$\wtt w(r)=2\log r+c+o(r). \tag 5.4$$

\noo (iii) In the case (5.4), by a blow-up argument and the
comparison principle, we have
$$w(x)=2\log |x|+c+o(1). \tag 5.5$$
Furthermore $w$ has isolated singularities.

\noo (iv) By (5.5) we prove that $w$ has at most one singular
point. Indeed, if $g\in [g_0]_\Ga$, the Ricci curvature $Ric_g\ge
0$, as indicated before Theorem 2.4. Hence the ratio $Q(r)=\frac
{Vol(B_{y, r}[g])}{Vol(B^e_{y, r})}$ is decreasing, where $B_{y,
r}[g]$ denotes geodesic ball in the metric $g$ and $B_{y, r}^e$
denotes ball in the Euclidean space. Hence
$$Q(r) \le 1. $$
By (5.5), we have
$$Q(r)\to m$$
as $r\to\infty$, where $m$ is the number of singular points. Hence
$m=1$ and
$$Q(r)\equiv 1.\tag 5.6$$

\noo (v) By (5.6) we show that if $w$ has a singular point
$\{0\}$, then $w\in C^\infty(\M\back\{0\})$. Indeed, for any point
$y\in\M\back\{0\}$, consider the Perron lifting of $w$ in $B_{y,
r}$, that is the function $w^*$ given by
$$\align
& \sigma_k(\lam(w^*))=0\ \ \text{in}\ \ B_{y, r}\\
& w^*=w\ \ \ \text{in}\ \ \M\back\{ B_{y, r}\}\\
\endalign $$
Then $w^*$ is admissible and $w^*\ge w$. By (5.6) we conclude
$w=w^*$ and so $w\in C^{1,1}$. By (5.6) we can show furthermore
that the scalar curvature $R_g=0$. Hence $w\in C^\infty
(\M\back\{0\})$.

\noo (vi) It follows that $(\M\back\{0\},  g)$ is a complete
smooth manifold with $Ric_g\ge 0$ and $Q(r)\equiv 1$. Hence it is
isometric to $\R^n$, which in turn implies that $\M$ is conformal
to $S^n$, a contradiction.

Having established the estimate (5.3), one can prove the existence
of solutions of (1.11) by a degree argument. We remark that under
the conditions in Theorem 5.1, there is a solution to
$$F(\D^2_{conf}v +\frac {n-2}{2}vA_{g_0}))=\phi v^p\tag 5.7$$
for any $p>\alpha\frac{n+2}{n-2}$. Hence there is no critical
exponent in this case.

\noo{\it Remark 5.1}. As in [TW2], one can prove that if
$\Ga\subset \Sigma_\delta$ for some $0\le \delta<\frac {1}{n-2}$,
and if $g=v^{\frac {4}{n-2}}g_0$ is an admissible metric on
$\R^n$, then either $v$ is Holder continuous or
$$v(x)=C|x|^{2-n}.\tag 5.8$$

\noo{\it Remark 5.2}. Our assumption on $f$ is stronger than that
in [L], except that we drop the smoothness condition $f\in C^1$ in
[L], so that it embraces example (ii). In the proof of the
Liouville theorem [L], the smoothness $f\in C^1$ is not required,
the local Lipschitz continuity of $f$, which follows from the
concavity of $f$, is enough. Alternatively, one can also use
smooth functions to approximate $f$, as the estimates (5.1) and
(5.3) are independent of $f$.  We also remark that if $v$ is an
admissible solution with bounded second derivatives, equation
(1.11) is uniformly elliptic and so $v\in C^{2, \alpha}$ [GT].

\vskip20pt

\baselineskip=12.0pt
\parskip=2.0pt

\Refs\widestnumber\key{ABC}

\item {[A]} T. Aubin,
       Some nonlinear problems in Riemannian geometry,
       Springer, 1998.

\item {[CNS]} L.A. Caffarelli, L. Nirenberg, and J. Spruck,
       Dirichlet problem for nonlinear second order
       elliptic equations III.
       Functions of the eigenvalues of the Hessian,
       Acta Math. 155(1985), 261--301.

\item {[CGY]} A. Chang, M. Gursky, P. Yang,
       An a priori estimate for a fully nonlinear equation on
       four-manifolds,
       J. Anal. Math. 87 (2002), 151--186.

\item {[Ch]} S. Chen,
        Local estimates for some fully nonlinear elliptic equation,
        preprint of September 2005.

\item {[CW]} K.S. Chou and X-J. Wang,
       A variational theory of the Hessian equation,
       Comm. Pure Appl. Math. 54 (2001), 1029--1064.

\item {[G]} C. Gerhardt,
       Closed Weingarten hypersurfaces in Riemannian manifolds,
       J. Diff. Geom., 43(1996), 612-641.

\item {[GT]}  D. Gilbarg and N.S. Trudinger,
       Elliptic partial differential equations of second order,
       Springer, 1983.

\item {[GW1]} P. Guan and G. Wang,
       Local estimates for a class of fully nonlinear equations
       arising from conformal geometry,
       Int. Math. Res. Not. (2003), 1413--1432.

\item {[GW2]} P. Guan and G. Wang,
        A fully nonlinear conformal flow on locally conformally
        flat manifolds,
        J. Reine Angew. Math. 557 (2003), 219--238.

\item {[GW3]} P. Guan and Xu-Jia Wang,
       On a Monge-Amp\`ere equation arising in geometric optics,
       J. Diff. Geom., 48(1998), 205--223.

\item {[GV1]} M. Gursky and J.Viaclovsky,
       Prescribing symmetric functions of the eigenvalues of the
       Ricci tensor, Ann. Math., to appear.

\item {[GV2]} M. Gursky and J.Viaclovsky,
       Convexity and singularities of curvature equations in
       conformal geometry, arXiv:math.DG/0504066.

\item {[I]} N. Ivochkina,
           Solution of the Dirichlet problem for some
           equations of Monge-Amp\`ere type.
           Mat. Sb.,  128(1985), 403-415.

\item {[Kr]}  N.V. Krylov,
       On the general notion of fully nonlinear second-order elliptic
       equations, Trans. Amer. Math. Soc., 347(1995), 857-895.

\item {[LL1]} A. Li and Y.Y. Li,
       On some conformally invariant fully nonlinear equations.
       Comm. Pure Appl. Math. 56 (2003), 1416--1464.

\item {[LL2]} A. Li and Y.Y. Li,
       On some conformally invariant fully nonlinear equations
       I\!I, Liouville, Harnack, and Yamabe, Acta Math., to
       appear.

\item {[L]} Y.Y. Li,
       Degenerate conformal invariant fully nonlinear elliptic
       equations, arXiv:math.AP/0504598.

\item {[STW]} W.M. Sheng, N.S. Trudinger, X.-J. Wang,
       The Yamabe problem for higher order curvatures,\newline
       arXiv:math.DG/0505463.

\item {[SY]} R. Schoen and S.T. Yau,
       Lectures on Differential geometry.
       International Press, 1994.

\item {[TW1]} N.S. Trudinger and X-J. Wang,
       Hessian measures I\!I,
       Ann. of Math. (2) 150 (1999), 579--604.

\item {[TW2]} N.S. Trudinger and X-J. Wang,
      On Harnack inequalities and singularities of
      admissible metrics in the Yamabe problem,
      arXiv:math.DG/0509341.

\item {[V]} J. Viaclovsky,
       Estimates and existence results for some fully nonlinear
       elliptic equations on Riemannian manifolds,
       Comm. Anal. Geom. 10 (2002), 815--846.

\item {[Ye]} R. Ye,
       Global existence and convergence of Yamabe flow.
       J. Differential Geom. 39 (1994),  35--50.

\endRefs

\enddocument
\end